\documentclass{article}
\usepackage{amsmath}

\input{tcilatex}
\begin{document}

\title{Another proof of Pell identities by using the determinant of
tridiagonal matrix}
\author{Meral Ya\c{s}ar\thanks{%
e-mail:myasar@nigde.edu.tr} \&Durmu\c{s} Bozkurt\thanks{%
e-mail:dbozkurt@selcuk.edu.tr} \and Department of Mathematics, Nigde
University and \and Department of Mathematics, Selcuk University}
\maketitle

\begin{abstract}
In this paper, another proof of Pell identities is presented by using the \
determinant of tridiagonal matrices. It is calculated via the Laplace
expansion.
\end{abstract}


\section{Introduction}

Pell numbers are defined as%
\begin{equation*}
P_{n}=2P_{n-1}+P_{n-2}
\end{equation*}%
with the initial conditions $P_{0}=0,P_{1}=1$ for $n\geq 2.$

In \cite{1}, a complex factorization formula for $(n+1)th$ Pell number is
obtained to provided the tridiagonal matrix%
\begin{equation*}
N(n)=\left[ 
\begin{array}{cccccc}
2i & 1 &  &  &  &  \\ 
1 & 2i & 1 &  &  &  \\ 
& 1 & 2i & 1 &  &  \\ 
&  & \ddots & \ddots & \ddots &  \\ 
&  &  & 1 & 2i & 1 \\ 
&  &  &  & 1 & 2i%
\end{array}%
\right]
\end{equation*}%
as in the following:%
\begin{equation}
P_{n+1}=m\left\vert N(n)\right\vert ,\ \ \ m=\left\{ 
\begin{array}{rc}
1, & n\equiv 0\ (\func{mod}4) \\ 
-i, & n\equiv 1\ (\func{mod}4) \\ 
-1, & n\equiv 2\ (\func{mod}4) \\ 
i, & n\equiv 3\ (\func{mod}4)%
\end{array}%
\right. .  \label{g}
\end{equation}%
In \cite{2}, an identity of Fibonacci numbers is proved via the determinant
of\linebreak the tridiagonal matrix. In \cite{3}, the authors showed the
connection between Fibonacci numbers and Chebyshev polynomials and obtained
a complex factorization for Fibonacci numbers by using a sequence of the
tridiagonal matrices. Then with a small difference in the tridiagonal
matrix, it is showed that how Lucas numbers and Chebyshev polynomials are
connected to each other. Two complex factorization are obtained by using the 
$n\times n$ tridiagonal and anti-tridiagonal matrix for $n$ is even in \cite%
{4}.

In this paper, we give another proof of Pell identities%
\begin{equation}
P_{2n}=P_{n}(P_{n+1}+P_{n-1})  \label{r}
\end{equation}%
\begin{equation}
P_{n}=P_{k}P_{n-k+1}+P_{k-1}P_{n-k}  \label{m}
\end{equation}%
where $k$ is a positive integer. For $1\leq k\leq n$%
\begin{equation}
\begin{array}{c}
P_{n}=P_{1}P_{n}+P_{0}P_{n-1} \\ 
P_{n}=P_{2}P_{n-1}+P_{1}P_{n-2} \\ 
P_{n}=P_{3}P_{n-2}+P_{2}P_{n-3} \\ 
P_{n}=P_{4}P_{n-3}+P_{3}P_{n-4} \\ 
\vdots \\ 
P_{n}=P_{n}P_{1}+P_{n-1}P_{0}.%
\end{array}
\label{n}
\end{equation}

\section{Main Result}

Let $A$ be an $n\times n$ matrix, $A([i_{1},i_{2},\ldots
,i_{k}],[j_{1},j_{2},\ldots ,j_{k}])$ be the $k\times k$ submatrix of $A$
and $M([i_{1},i_{2},\ldots ,i_{k}],[j_{1},j_{2},\ldots ,j_{k}])$ be the $%
(n-k)\times (n-k)$ minor of the matrix $A.$ The cofactor of $A$ is defined by%
\begin{equation*}
\mathring{A}([i_{1},i_{2},\ldots ,i_{k}],[j_{1},j_{2},\ldots
,j_{k}])=(-1)^{m}M([i_{1},i_{2},\ldots ,i_{k}],[j_{1},j_{2},\ldots ,j_{k}])
\end{equation*}%
where $1\leq i_{1},i_{2},\cdots ,i_{k}\leq n$ and $m=\tsum%
\nolimits_{r=1}^{k}(i_{r}+j_{r}).$

The determinant of the matrix $A$ is%
\begin{eqnarray*}
\det (A) &=&\dsum\limits_{1\leq i_{1},i_{2},\cdots ,i_{k}\leq n}\det
(A([i_{1},i_{2},\ldots ,i_{k}],[j_{1},j_{2},\ldots ,j_{k}])) \\
&&\times \det (\mathring{A}([i_{1},i_{2},\ldots ,i_{k}],[j_{1},j_{2},\ldots
,j_{k}])).
\end{eqnarray*}%
If $A(i,j)=a_{ij},$ then $\mathring{A}(i,j)=(-1)^{i+j}M(i,j)=\mathring{A}%
_{ij}$ and the determinant is%
\begin{equation*}
\det (A)=\dsum\limits_{i=1}^{n}a_{ij}\mathring{A}_{ij}.
\end{equation*}%
This is the famous Laplace expansion formula \cite{2}. We will use this
formula to proof of Pell identities in (\ref{m}).

The cofactors of the first row of the matrix $N(n)$ are 
\begin{equation*}
\mathring{A}_{11}=\left\{ 
\begin{array}{rc}
-P_{n-1} & ,n\equiv 0(\func{mod}4) \\ 
-iP_{n-1} & ,n\equiv 1(\func{mod}4) \\ 
P_{n-1} & ,n\equiv 2(\func{mod}4) \\ 
iP_{n-1} & ,n\equiv 3(\func{mod}4)%
\end{array}%
\right.
\end{equation*}%
\begin{equation*}
\mathring{A}_{12}=\left\{ 
\begin{array}{rc}
-iP_{n-2} & ,n\equiv 0(\func{mod}4) \\ 
P_{n-2} & ,n\equiv 1(\func{mod}4) \\ 
iP_{n-2} & ,n\equiv 2(\func{mod}4) \\ 
-P_{n-2} & ,n\equiv 3(\func{mod}4)%
\end{array}%
\right. .
\end{equation*}%
By using the Laplace expansion formula the determinant of the matrix $N(n-1)$
is%
\begin{equation}
\det (N(n-1))=2P_{n-1}+P_{n-2}.  \label{a}
\end{equation}

From (\ref{g}), the initial value $P_{0}=0,P_{1}=1,P_{2}=2$ and the fact $%
P_{2}=2P_{1}$ are used in (\ref{a}), then we have%
\begin{eqnarray}
P_{n} &=&2P_{1}P_{n-1}+P_{1}P_{n-2}  \notag \\
&=&P_{2}P_{n-1}+P_{1}P_{n-2}.  \label{c}
\end{eqnarray}

If the first two rows of the matrix $N(n-1)$ are chosen, there are only
three $2\times 2$ submatrices of the matrix $N(n-1)$ whose determinants are
nonzero. i.e.%
\begin{eqnarray*}
A([1,2],[1,2]) &=&\left[ 
\begin{array}{cc}
2i & 1 \\ 
1 & 2i%
\end{array}%
\right] =-P_{3} \\
A([1,2],[1,3]) &=&\left[ 
\begin{array}{cc}
2i & 0 \\ 
1 & 1%
\end{array}%
\right] =iP_{2} \\
A([1,2],[2,3]) &=&\left[ 
\begin{array}{cc}
1 & 0 \\ 
2i & 1%
\end{array}%
\right] =P_{1}
\end{eqnarray*}%
and their cofactors are%
\begin{equation*}
\mathring{A}([1,2],[1,2])=\left\{ 
\begin{array}{rc}
iP_{n-2} & ,n\equiv 0(\func{mod}4) \\ 
-P_{n-2} & ,n\equiv 1(\func{mod}4) \\ 
-iP_{n-2} & ,n\equiv 2(\func{mod}4) \\ 
P_{n-2} & ,n\equiv 3(\func{mod}4)%
\end{array}%
\right.
\end{equation*}%
\begin{equation*}
\mathring{A}([1,2],[1,3])=\left\{ 
\begin{array}{rc}
-P_{n-3} & ,n\equiv 0(\func{mod}4) \\ 
-iP_{n-3} & ,n\equiv 1(\func{mod}4) \\ 
P_{n-3} & ,n\equiv 2(\func{mod}4) \\ 
iP_{n-3} & ,n\equiv 3(\func{mod}4)%
\end{array}%
\right.
\end{equation*}%
\begin{equation*}
\mathring{A}([1,2],[2,3])=0.
\end{equation*}%
By using the Laplace expansion the determinant of the matrix $N(n-1)$ is%
\begin{equation}
\det (N(n-1))=\left\{ 
\begin{array}{rc}
-iP_{3}P_{n-2}-iP_{2}P_{n-3} & ,n\equiv 0(\func{mod}4) \\ 
P_{3}P_{n-2}+P_{2}P_{n-3} & ,n\equiv 1(\func{mod}4) \\ 
iP_{3}P_{n-2}+iP_{2}P_{n-3} & ,n\equiv 2(\func{mod}4) \\ 
-P_{3}P_{n-2}-P_{2}P_{n-3} & ,n\equiv 3(\func{mod}4)%
\end{array}%
\right. .  \label{b}
\end{equation}

From (\ref{g}) and(\ref{b}), we obtain%
\begin{equation}
P_{n}=P_{3}P_{n-2}+P_{2}P_{n-3}.  \label{d}
\end{equation}

If the first three rows of the matrix $N(n-1)$ are chosen, there are only
four $3\times 3$ submatrices of the matrix $N(n-1)$ whose determinants are
nonzero:%
\begin{eqnarray*}
A([1,2,3],[1,2,3]) &=&\left[ 
\begin{array}{ccc}
2i & 1 & 0 \\ 
1 & 2i & 1 \\ 
0 & 1 & 2i%
\end{array}%
\right] =-iP_{4} \\
A([1,2,3],[1,2,4]) &=&\left[ 
\begin{array}{ccc}
2i & 1 & 0 \\ 
1 & 2i & 0 \\ 
0 & 1 & 1%
\end{array}%
\right] =-P_{3} \\
A([1,2,3],[1,3,4]) &=&\left[ 
\begin{array}{ccc}
2i & 0 & 0 \\ 
1 & 1 & 0 \\ 
0 & 2i & 1%
\end{array}%
\right] =iP_{2} \\
A([1,2,3],[2,3,4]) &=&\left[ 
\begin{array}{ccc}
1 & 0 & 0 \\ 
2i & 1 & 0 \\ 
1 & 2i & 1%
\end{array}%
\right] =P_{1}
\end{eqnarray*}%
and their cofactors are%
\begin{equation*}
\mathring{A}([1,2,3],[1,2,3])=\left\{ 
\begin{array}{rc}
P_{n-3} & ,n\equiv 0(\func{mod}4) \\ 
iP_{n-3} & ,n\equiv 1(\func{mod}4) \\ 
-P_{n-3} & ,n\equiv 2(\func{mod}4) \\ 
-iP_{n-3} & ,n\equiv 3(\func{mod}4)%
\end{array}%
\right.
\end{equation*}%
\begin{equation*}
\mathring{A}([1,2,3],[1,2,4])=\left\{ 
\begin{array}{rc}
iP_{n-4} & ,n\equiv 0(\func{mod}4) \\ 
-P_{n-4} & ,n\equiv 1(\func{mod}4) \\ 
-iP_{n-4} & ,n\equiv 2(\func{mod}4) \\ 
P_{n-4} & ,n\equiv 3(\func{mod}4)%
\end{array}%
\right.
\end{equation*}%
\begin{equation*}
\mathring{A}([1,2,3],[1,3,4])=0
\end{equation*}%
\begin{equation*}
\mathring{A}([1,2,3],[2,3,4])=0
\end{equation*}%
By using the Laplace expansion the determinant of the matrix $N(n-1)$ is
evaluated as follows:%
\begin{equation}
\det (N(n-1))=\left[ 
\begin{array}{rc}
-iP_{4}P_{n-3}-iP_{3}P_{n-4} & ,n\equiv 0(\func{mod}4) \\ 
P_{4}P_{n-3}+P_{3}P_{n-4} & ,n\equiv 1(\func{mod}4) \\ 
iP_{4}P_{n-3}+iP_{3}P_{n-4} & ,n\equiv 2(\func{mod}4) \\ 
-P_{4}P_{n-3}-P_{3}P_{n-4} & ,n\equiv 3(\func{mod}4)%
\end{array}%
\right. .  \label{e}
\end{equation}

From (\ref{g}) and(\ref{e}), we have%
\begin{equation}
P_{n}=P_{4}P_{n-3}+P_{3}P_{n-4}.  \label{f}
\end{equation}%
The remainig identities in (\ref{n}) can be shown similarly.

Now, we give another proof of following Pell identity:%
\begin{equation}
P_{2n}=P_{n}(P_{n+1}+P_{n-1}).  \label{h}
\end{equation}%
If we choose the first $(n-1)$ rows of the matrix $N(2n-1),$ there are only $%
n$ the $(n-1)\times (n-1)$ submatrices of the matrix $N(2n-1)$ whose
determinants are nonzero but only the cofactors of two of them are nonzero.
i.e.%
\begin{equation*}
A([1,2,\ldots ,n-1],[1,2,\ldots ,n-2,n-1])=\left\{ 
\begin{array}{rc}
-iP_{n} & ,n\equiv 0(\func{mod}4) \\ 
P_{n} & ,n\equiv 1(\func{mod}4) \\ 
iP_{n} & ,n\equiv 2(\func{mod}4) \\ 
-P_{n} & ,n\equiv 3(\func{mod}4)%
\end{array}%
\right.
\end{equation*}%
\begin{equation*}
A([1,2,\ldots ,n-1],[1,2,\ldots ,n-2,n])=\left\{ 
\begin{array}{rc}
-P_{n-1} & ,n\equiv 0(\func{mod}4) \\ 
-iP_{n-1} & ,n\equiv 1(\func{mod}4) \\ 
P_{n-1} & ,n\equiv 2(\func{mod}4) \\ 
iP_{n-1} & ,n\equiv 3(\func{mod}4)%
\end{array}%
\right.
\end{equation*}%
and their cofactors are%
\begin{equation*}
\mathring{A}([1,2,\ldots ,n-1],[1,2,\ldots ,n-2,n-1])=\left\{ 
\begin{array}{rc}
P_{n+1} & ,n\equiv 0(\func{mod}4) \\ 
iP_{n+1} & ,n\equiv 1(\func{mod}4) \\ 
-P_{n+1} & ,n\equiv 2(\func{mod}4) \\ 
-iP_{n+1} & ,n\equiv 3(\func{mod}4)%
\end{array}%
\right.
\end{equation*}%
\begin{equation*}
\mathring{A}([1,2,\ldots ,n-1],[1,2,\ldots ,n-2,n])=\left\{ 
\begin{array}{rc}
iP_{n} & ,n\equiv 0(\func{mod}4) \\ 
-P_{n} & ,n\equiv 1(\func{mod}4) \\ 
-iP_{n} & ,n\equiv 2(\func{mod}4) \\ 
P_{n} & ,n\equiv 3(\func{mod}4)%
\end{array}%
\right. .
\end{equation*}%
From Laplace expansion the determinant of the matrix $N(2n-1)$ is%
\begin{equation}
\det (N(2n-1))=\left\{ 
\begin{array}{rc}
-iP_{n}(P_{n+1}+P_{n-1}) & ,n\equiv 0,2(\func{mod}4) \\ 
iP_{n}(P_{n+1}+P_{n-1}) & ,n\equiv 1,3(\func{mod}4)%
\end{array}%
\right. .  \label{s}
\end{equation}%
From (\ref{g}) and (\ref{s}) we have%
\begin{equation*}
P_{2n}=P_{n}(P_{n+1}+P_{n-1}).
\end{equation*}%
Thus the proof is completed.

\end{document}